\documentclass[12pt]{article}

\usepackage{amssymb, amsmath, amsthm, amsfonts}
\usepackage{multirow}
\usepackage{booktabs}
\usepackage{subfigure}
\usepackage[authoryear]{natbib}
\usepackage[colorlinks,citecolor=blue,urlcolor=blue]{hyperref}
\usepackage{graphicx}
\usepackage{newtxtext}
\usepackage{epsfig,latexsym,verbatim}
\usepackage{bm}
\usepackage{mathrsfs}
\usepackage{multicol}
\usepackage{enumitem}
\usepackage{makecell}
\usepackage[figuresright]{rotating}
\usepackage{xr}

\newcommand{\blind}{0}

\theoremstyle{plain}

\newtheorem{theorem}{Theorem}
\newtheorem{example}{Example}
\newtheorem{proposition}{Proposition}
\newtheorem{lemma}{Lemma}[section]
\newtheorem{remark}{Remark}
\newtheorem{corollary}{Corollary}
\theoremstyle{remark}

\def\P{\mathbb{P}}
\def\E{\mathbb{E}}
\def\R{\mathbb{R}}
\def\Var{\mathrm{Var}}
\def\Cov{\mathrm{Cov}}
\def\o{{\scriptstyle{\mathcal{O}}}}
\def\O{\mathcal{O}}

\renewcommand{\footnote}{}

\DeclareMathOperator{\Tr}{Tr}

\begin{document}

\def\spacingset#1{\renewcommand{\baselinestretch}%
{#1}\small\normalsize} \spacingset{1}

\date{}
\if0\blind
{
  \title{\bf Strong Gaussian approximation for U-statistics in high dimensions and beyond}
  \author{Weijia Li\footnotemark[1] \footnotemark[3],$\,$ Leheng Cai\footnotemark[1] \footnotemark[3], $\,$ and$\,$ Qirui Hu\footnotemark[2] \footnotemark[4] \hspace{.2cm}}
  \maketitle
  \renewcommand{\thefootnote}{\fnsymbol{footnote}}

  \footnotetext[1]{Tsinghua University}
  \footnotetext[2]{Shanghai University of Finance and Economics}
  \footnotetext[3]{Co-first authors. These authors contributed equally to this work.}
  \footnotetext[4]{Corresponding author: huqirui@mail.shufe.edu.cn}
}

\if1\blind
{
  \bigskip
  \bigskip
  \bigskip
  \begin{center}
    { \bf Strong Gaussian approximation for U-statistics in high dimensions and beyond}
  \end{center}
  \medskip
} \fi

\bigskip

\abstract{We establish a strong Gaussian approximation for high-dimensional non-degenerate U-statistics with diverging dimension. Under mild assumptions, we construct, on a sufficiently rich probability space, a Gaussian process that uniformly approximates the entire sequential U-statistic process. The approximation error is explicitly characterized and vanishes under polynomial growth of the dimension. The key technical contribution is a sharp martingale maximal inequality for completely degenerate U-statistics, combined with a high-dimensional strong approximation for independent sums. This coupling yields functional Gaussian limits without relying on $\mathcal{L}^\infty$-type bounds or bootstrap arguments. The theory is illustrated through three representative examples of U-statistics: the spatial Kendall’s tau matrix, the multivariate Gini’s mean difference, and the characteristic dispersion parameter. As applications, we derive Brownian bridge approximations for U-statistic-based change-point statistics and develop a self-normalized relevant testing procedure whose limiting distribution is fully pivotal. The framework naturally accommodates bounded kernels and therefore remains valid under heavy-tailed distributions. Overall, our results provide a unified probability-theoretic foundation for high-dimensional inference based on U-statistics.}

\vspace{0.5cm}
\noindent\textbf{Keywords: change point test, high-dimensional statistics, relevant hypotheses, strong Gaussian approximation, U-statistics}

\newpage
\spacingset{1.9}

\section{Introduction}\label{Sec:intro}
$U$-statistics, introduced by \citet{Hoeffding1948}, provide unbiased estimators for a broad class of parameters of the form $\bm{\theta}=\E[\bm{h}(X_1,X_2)]$, where $\bm{h}$ is a symmetric kernel and $X_1,X_2$ are independent copies of an observation. In contemporary applications, $\bm{h}$ is frequently vector valued and the dimension $d$ of the target parameter grows with the sample size. Representative examples include robust dependence and dispersion measures such as spatial Kendall's tau matrices \citep{Visuri} and rank- or sign-based high-dimensional procedures used in robust inference \citep{hanfangPCA}. In high-dimensional regimes, pairwise kernels are particularly attractive because they can be chosen to be bounded or Lipschitz, leading to procedures that remain well defined even when the underlying distribution is heavy tailed.

A recurring methodological step in these problems is a Gaussian approximation of the centered and scaled $U$-statistic. Beyond weak convergence, it is often advantageous to construct an   coupling between the statistic and a Gaussian counterpart on the same  probability space. Such a strong Gaussian approximation (also called a strong invariance principle) yields distributional approximations for a large class of functionals and is indispensable in sequential problems such as change-point analysis and self-normalized inference. In this paper we study strong Gaussian approximation for the sequential process of high-dimensional $U$-statistics in the Euclidean norm. Let $\bm{U}_k$ denote the $U$-statistic computed from the first $k$ observations and let $\bm{T}_k$ be the associated $\sqrt{n}$-scaled and centered process (defined in Section~\ref{sec:sequential}). Our goal is to construct a Gaussian partial-sum process  $\bm{W}_k$ in a richer probability space, such that
\[
\max_{2\le k\le n}\|\bm{T}_k-\bm{W}_k\|_2
\]
is asymptotically ignorable. By Strassen's theorem (see, e.g., \citet{Pollard_2002}), such a coupling control 
yields a uniform bound on the L\'evy–Prohorov distance 
between the corresponding laws, which provides an intuitive perspective on Gaussian coupling from the viewpoint of probabilistic approximation. 

Strong approximation has a long history for sums of independent random variables. The Hungarian construction of \citet{KomlosMajorTusnady1975,KomlosMajorTusnady1976} yields nearly optimal almost sure approximations in one dimension, with multidimensional extensions and approximation theorems for random vectors developed in, among others, \citet{BerkesPhilipp1979,Yurinskii1977,Zaitsev1998}. These results and their refinements concern (possibly vector-valued) partial sums, and they constitute the classical foundation for strong invariance principles. More recently, \citet{mies2023} established a sequential strong Gaussian approximation for partial sums of high-dimensional random vectors in Euclidean norm with explicit rates, which will serve as a key ingredient in our treatment of the linear (H\'ajek projection) component of the $U$-statistic.

For $U$-statistics, functional limit theorems and invariance principles in fixed dimension have been studied extensively. The nondegenerate case is treated in \citet{Hall1979}, while degenerate kernels lead to different limiting objects as in \citet{Neuhaus1977}. Almost sure behavior and strong approximations were obtained in early work such as \citet{Sen1974} and in the Kiefer-process based approximation of \citet{DehlingDenkerPhilipp1984}. While these classical works provide deep insights into the pathwise behavior of $U$-statistics, they are largely developed for fixed-dimensional targets and do not address the high-dimensional regime $d=d_n\to\infty$ where explicit approximation rates and uniform-in-time control become essential.

In high-dimensional statistics, a different line of research has developed nonasymptotic distributional Gaussian approximations, primarily for maxima and hyperrectangles. Starting from the high-dimensional CLT and multiplier bootstrap theory for sums of random vectors \citep{ChernozhukovChetverikovKato2013,ChernozhukovChetverikovKato2017}, analogous results have been obtained for $U$-statistics and related objects; see, for example, \citet{chen-Ustat,ChenKato2019,SongChenKato2019,ChengLiuPeng2022,ImaiKoike2025}. These CCK-type bounds are powerful for simultaneous inference in the $\mathcal{L}^\infty$-geometry and can allow $d$ to be very large, but they are geared towards distributional approximation of maximum-type functionals. Translating such results into a sequential strong coupling in Euclidean norm is nontrivial and is not covered by the existing literature. In particular, our focus on the $\mathcal{L}^2$-geometry naturally leads to polynomial dimension growth and calls for different techniques than anti-concentration and $\mathcal{L}^\infty$-based empirical process bounds.

Our main contribution is to establish a sequential strong Gaussian approximation for high-dimensional, vector-valued $U$-statistics of order two. Under a finite $q>2$ moment condition on the H\'ajek projection and only a finite second moment for the degenerate kernel, we construct on an enriched probability space independent Gaussian vectors $\{\bm Z_i\}_{i=1}^n$ with matching covariances such that the associated Gaussian partial sums $\bm W_k=n^{-1/2}\sum_{i=1}^k\bm Z_i$ satisfy
\[
\max_{2\le k\le n}\|\bm{T}_k-\bm{W}_k\|_2
= \O_p\!\left(B\sqrt{\log n}\left( {d}/{n}\right)^{ 1/4- {1}/{(2q)}}\right),
\]
and the right-hand side vanishes when $d$ grows at an explicit polynomial rate (cf.\ Assumption~(A2) in Section~\ref{sec:sequential}). We also provide a Gaussian approximation for the global statistic in an independent but not necessarily identically distributed setting.

A key byproduct is a maximal inequality for vector-valued degenerate $U$-statistics showing that the degenerate remainder is uniformly of order $\sqrt{d\log n}$ after normalization. This step is technically delicate because the degenerate component is neither a sum of independent terms nor a standard empirical process indexed by a fixed class, and the sequential setting requires uniform control over $k=2,\ldots,n$. We overcome this difficulty by embedding the degenerate sequential $U$-statistic into a martingale with respect to the natural filtration and then applying a vector-valued martingale maximal inequality \citep{Baiinequality} together with a classical martingale inequality \citep{chow1960martingale}. This martingale route avoids higher-order moment and tail assumptions and is crucial for robust kernels; in particular, bounded kernels imply the required moments regardless of the tail behavior of $X_i$.

We demonstrate the scope of the strong approximation through two applications. In Section~\ref{sec:relevant}, we develop a self-normalized test for relevant hypotheses for high-dimensional parameter $\bm{\theta}$ with the form $\|\bm{\theta}-\bm{\theta}_0\|_2^2\le\Delta$, obtaining an asymptotically pivotal limit without estimating a high-dimensional covariance matrix.  In Section~\ref{sec:changepoint} we study retrospective change-point detection based on $U$-statistic CUSUM processes and derive a Brownian-bridge limit under the null hypothesis, together with consistency under alternatives. These results complement recent work on $U$-statistic based change-point methods such as \citet{wegner2023robust} and the most recent results in  \citet{dehling2026}, by providing a sequential strong approximation in Euclidean norm with explicit dimension dependence. 

The rest of the paper is organized as follows. Section~\ref{sec:sequential} presents preliminaries and our main results for strong Gaussian approximation. Section~\ref{Sec:statisticsapplication} develops the statistical applications with established theories. All detailed proofs are are collected in supplemental materials.

\section{Main results}\label{sec:sequential}
 
\subsection{Preliminaries}\label{sec:main-pre}

Let $\{X_i\}_{i=1}^n$ be i.i.d.\ random vectors in $\R^p$ with distribution $F$. Let $\bm{h}=(h_1,\dots,h_d)^\top:\R^p\times\R^p\to\R^d$ be a measurable symmetric kernel, i.e.\ $\bm{h}(x_1,x_2)=\bm{h}(x_2,x_1)$.   The parameter we are interested in is $\bm \theta := \E[\bm{h}(X_1,X_2)]\in\R^d$.
The classical multivariate U-statistic is defined as
\[
\bm{U}_n := \frac{2}{n(n-1)}\sum_{1\le i<j\le n} \bm{h}(X_i,X_j),\quad n\ge 2.
\]
By the Hoeffding decomposition, 
\begin{align}
    \bm{U}_n-\bm{\theta}=\frac{2}{n}\sum_{i=1}^n \bm{g}(X_i)+\frac{1}{n(n-1)}\sum_{1\le i\neq j\le n} \bm{f}(X_i,X_j),
    \label{eq:Hoeffding-deco}
\end{align}
where $\bm{g}(X)=\E\left\{\bm{h}(X,X')|X\right\}-\bm \theta$ and $\bm{f}(X,X')=\bm{h}(X,X')-\bm{g}(X)-\bm{g}(X')-\bm\theta$. Here,  $X$ and $X'$ are independent and have the same distribution as $X_i$. 
Moreover, $\bm{g}(\cdot)$ and $\bm{f}(\cdot,\cdot)$ are $d$-dimensional functions consisting of $(g_1,\ldots,g_d)$ and $(f_1,\ldots,f_d)$, respectively.

An advantage of U-statistics is   their ability to capture complex nonlinear pairwise interactions between observations through kernel functions. To further illustrate the multivariate U-statistic and our main results, we present several examples below. 
\begin{example}[Multivariate Gini's Mean Difference]\label{exa:GMD}
In various disciplines such as economics and social sciences, Gini's Mean Difference (GMD) serves as a fundamental measure of statistical dispersion. Historically utilized to quantify  wealth inequality within a population (e.g., \cite{yitzhaki2013gini}), the GMD conceptually captures the expected absolute difference between any two independent individuals. 

Compared with the classical sample variance, the  GMD  provides a more robust measure of dispersion. While the sample variance relies on squared deviations and is therefore highly sensitive to extreme observations, the GMD is based on pairwise absolute differences, which makes it  particularly appealing for heavy-tailed distributions. 
Formally, one defines high-dimensional target parameter component-wise as
\[
\theta_\ell := \E|X_{1\ell} - X_{2\ell}|, \quad \ell = 1, \dots, d.
\]
To estimate this parameter without bias, we naturally adopt the absolute difference kernel
$h_\ell(X_i, X_j) = |X_{i\ell} - X_{j\ell}|$. 
In the robust statistics literature, such U-statistics have been widely used to construct nonparametric and robust covariance analogues; see \cite{gini1} for instance. 

\end{example}

\begin{example}[Characteristic dispersion parameter]\label{exa:CDP}
While the Gini mean difference relaxes the moment condition to the first order, certain high-dimensional applications, such as high-frequency finance or robust signal processing, often involve extremely heavy-tailed data (e.g., Cauchy distributions), where even the first moment may not exist. To quantify multivariate dispersion in such situations, we adopt an approach based on the empirical characteristic function, motivated by the classical literature on stable laws (e.g., \cite{press1972}). 

For a random variable $X$, let $\varphi(t) = \E[\exp(\sqrt{-1}tX)]$ be its characteristic function. If $X_1, X_2$ are i.i.d.\ copies, the c.f. of their difference evaluated at $t=1$ yields a moment-free measure of dispersion.
$\E[\cos(X_1 - X_2)] = |\varphi(1)|^2$. 
One naturally defines the target parameter component-wise as:
$$ \theta_\ell := \E \cos(X_{i\ell} - X_{j\ell}) , \quad \text{for } \ell = 1, \dots, d. $$
This parameter     captures the generalized dispersion or scale of the distribution, reflecting the decay rate of the probability density in the characteristic domain. Notably, it is  strictly bounded in $(0, 1]$ for any probability distribution, completely bypassing the need for moment assumptions. For instance:
\begin{itemize}
    \item[(i)] If $X_{i\ell} \sim \mathcal{N}(\mu_\ell, \sigma_\ell^2)$, then $\theta_\ell = \exp(-\sigma_\ell^2)$, a monotone transformation of the variance under normality.
    \item[(ii)] If $X_{i\ell} \sim \text{Cauchy}(\mu_\ell, \gamma_\ell)$, then $\theta_\ell = \exp(-2\gamma_\ell)$, capturing the scale parameter $\gamma_\ell$.
    \item[(iii)] If $X_{i\ell} \sim \text{Laplace}(\mu_\ell, b_\ell)$, then $\theta_\ell = (1+b_\ell^2)^{-2}$.
\end{itemize} 
To estimate this parameter without bias, one adopts the   cosine kernel $h_\ell(X_i, X_j) = \cos(X_{i\ell} - X_{j\ell})$. 


\end{example}

\begin{example}[Spatial Kendall's tau matrix]\label{exa:SKT}
Constructing reliable gene co-expression networks from high-throughput biological data presents a unique statistical challenge. Unlike idealized multivariate settings, genomic readouts (e.g., RNA-seq or microarray expressions) are intrinsically affected by severe measurement errors, technical noise, and arbitrary scale variations. To capture the intrinsic spatial association while remaining invariant to such data contamination, biologists often favor the spatial Kendall's tau matrix, which is originally introduced by \cite{Visuri} and proved to be effective in robust principal component analysis like in \cite{hanfangPCA}. 
Unlike covariance-based measures, this parameter depends only on the signs of pairwise differences rather than their magnitudes. 

By vectorizing the matrix, the $d$-dimensional target parameter (where $d = p^2$) is   defined as:
$$ \bm{\theta} := \E\left[ \text{vec}\left( \frac{(X_1 - X_2)(X_1 - X_2)^\top}{\|X_1 - X_2\|_2^2} \right) \right]. $$
To estimate this parameter without bias, we employ the symmetric kernel:
$$ \bm{h}(X_i, X_j) = \text{vec}\left( \frac{(X_i - X_j)(X_i - X_j)^\top}{\|X_i - X_j\|_2^2} \right), \quad \text{for } i\neq j.$$
For any $i\neq j$, the unvectorized rank-one matrix $\frac{(X_i - X_j)(X_i - X_j)^\top}{\|X_i - X_j\|_2^2}$ has a trace of exactly $1$, which guarantees that the Euclidean norm of the vectorized kernel is strictly unity, i.e., $\|\bm{h}(X_i, X_j)\|_2 \equiv  1$. By mapping the pairwise differences onto a unit sphere, the magnitude of any measurement outlier is completely neutralized, retaining only the relative directional information of the genetic expression profiles. 

In fact, the spatial Kendall's tau matrix belongs to a general class of U-type robust covariance estimators, whose kernel function is formulated as
$$ h(X_i,X_j)=   \psi_{\xi}\left( \frac{1}{2} \|X_i - X_j\|_2^2 \right) \frac{(X_i - X_j)(X_i - X_j)^\top}{\|X_i - X_j\|_2^2}, $$
where $\psi_{\xi}(t) = \xi \psi(t/\xi)$, and $\psi(\cdot)$ is a robust score function (e.g., the Huber score), modulating the impact of extreme pairwise distances. The spatial Kendall's tau is a specific case of this family with $\psi_{\xi}(\cdot) \equiv 1$. For comprehensive theoretical properties of this general class, see \cite{YuRenZhou2024}.
\end{example}


\subsection{Sequential Gaussian approximation for partial sums}

While U-statistics based on a fixed sample size are fundamental, applications such as change point analysis discussed later in Section \ref{sec:changepoint} naturally require monitoring sequential data. Motivated by this, we investigate the asymptotic properties of the sequential partial sums of U-statistics in this work.

Let $\bm{U}_k$ be the U-statistic based on the first $k$ observations. To construct a sequential process that allows for asymptotic analysis, we define the scaled sequential statistic $\bm{T}_k$  in terms of the kernel $\bm{h}$: 
\begin{align}\label{eq:Tk_def}
    \bm{T}_k := \frac{k}{2\sqrt{n}}(\bm{U}_k - \bm{\theta}) = \frac{1}{\sqrt{n}(k-1)}\sum_{1\le i<j\le k} (\bm{h}(X_i,X_j)-\bm \theta) , \quad \text{for } 2\le k\le n.
\end{align}
Note that the conventional global statistic evaluated on the full sample $\bm{T}_n =  {\sqrt{n}} (\bm{U}_n - \bm{\theta})/2$ is   a special case of this sequential process at the terminal time point $n$.  

Our primary objective is to approximate the sequential process 
$\{\bm{T}_k\}_{k=2}^n$ by a sequence of Gaussian random vectors 
$\{\bm{W}_k\}_{k=2}^n$ through a uniform control of the physical distance 
$\|\bm{T}_k-\bm{W}_k\|_2$ on a possibly richer probability space. 
By Strassen's theorem \citep{Pollard_2002}, such a coupling yields a uniform bound on the Lévy–Prohorov distance 
$\rho(\bm{T}_k,\bm{W}_k,\epsilon)$ defined below:
\begin{align*}
    \rho(\bm{T}_k,\bm{W}_k,\epsilon)=\sup_{A\in \mathcal{B}}\max\left\{\P(\bm{T}_k\in A)-\P(\bm{W}_k\in A^{\epsilon}), \P(\bm{W}_k\in A)-\P(\bm{T}_k\in A^{\epsilon})\right\},
\end{align*}
where $A^{\epsilon}$ denotes the $\epsilon$-neighborhood of $A$, and $\mathcal B$ is  the collection of all Borel sets. 

To facilitate this bound, we apply the Hoeffding decomposition to separate the statistic into its linear projection and degenerate remainder:
\begin{align*}
    \bm{T}_k = \frac{1}{\sqrt{n}}\sum_{i=1}^k \bm{g}(X_i) + \frac{1}{2\sqrt{n}(k-1)}\sum_{1\le i \neq j\le k}\bm{f}(X_i,X_j),
\end{align*}
where the centered first-order projection $\bm{g}(\cdot)$ and the completely degenerate kernel $\bm{f}(\cdot,\cdot)$ are defined  as in \eqref{eq:Hoeffding-deco}. We denote the covariance matrix of the first order projection as $\bm{\Sigma}= \Cov\{\bm{g}(X_1)\}=\Cov\{\bm{h}(X_1,X_2),\bm{h}(X_1,X_3)\}$. The theoretical validity of this uniform strong approximation rests on the following regularity conditions:


\begin{itemize}
    \item[(A1)] There exists some $q>2$ such that $\left(\E\|\bm{g}(X)\|^q\right)^{1/q}\le B$.
    \item[(A2)] Let $D_r = \max_{1\le m \le d} (\E|f_{m}(X_1, X_2)|^r)^{1/r}$, and assume that $D_2 < \infty$.
    \item[(A3)] The bound $B$ in (A1) satisfies $B\asymp \sqrt{d}$ and the dimension $d=\o\left(n^{\frac{q-2}{3q-2}-\gamma}\right)$ for some small $\gamma>0$.
\end{itemize}

While Assumptions (A1) and (A2) are formulated on the Hoeffding projections $\bm{g}(\cdot)$ and $\bm{f}(\cdot,\cdot)$ to facilitate asymptotic derivations, they are   implied by mild moment conditions on the original kernel $\bm{h}(\cdot,\cdot)$ via standard moment inequalities. Specifically, a sufficient condition for (A1) is $(\E\|\bm{h}(X_1, X_2)\|^q)^{1/q} = \O(B)$, which typically holds when the components of $\bm{h}(\cdot,\cdot)$ have uniformly bounded $q$-th moments. Similarly, a sufficient condition for (A2) is  $\max_{1\le m \le d} \E|h_{m}(X_1, X_2)|^2 < \infty$. 
It is worth noting that in typical high-dimensional settings, the $q$-th moment bound $B$ in Assumption (A1) inherently scales with the dimension as $B \asymp \sqrt{d}$. 
\begin{theorem}\label{thm:sequential-borel}
    Under Assumptions (A1)-(A2),  there exists a sequence of independent Gaussian random vectors $\{\bm Z_i\}_{i=1}^n$ with $\bm Z_i \sim \mathcal{N}(\bm{0}, \bm{\Sigma})$ such that, for the Gaussian partial sum process $\bm{W}_k =  \sum_{i=1}^k \bm Z_i/\sqrt{n}$,  
    $$
    \mathbb{E} \left\{ \max_{2 \le k \le n} \left\| \bm{T}_k - \bm{W}_k \right\|_2 \right\} = \O\left( B \sqrt{\log n} \left(  {d}/{n} \right)^{  {1}/{4} -  {1}/{(2q)}} \right).
    $$
    Under (A3), the above approximation error vanishes asymptotically.
\end{theorem}
  Our sequential Gaussian approximation framework provides a methodological complement to the seminal max-type ($\mathcal{L}^\infty$-norm) Gaussian approximations pioneered by \cite{ChernozhukovChetverikovKato2013} for independent sums and later extended to U-statistics by \cite{chen-Ustat}. The main technical tools in this line of work involve Gaussian comparison and anti-concentration inequalities. While the max-type Gaussian approximation accommodates exponential dimension growth (e.g., $d = \mathcal{O}(\exp(n^c))$), it only bounds the Kolmogorov distance over hyperrectangles and therefore primarily yields a distributional approximation for the $\mathcal{L}^\infty$-norm of static statistics at a fixed sample size $n$. 
   In contrast, the proposed framework establishes a strong approximation uniformly over the entire sequence. Although the accumulation of approximation errors across all coordinates restricts the dimension to grow at a polynomial rate, this $\mathcal{L}^2$ formulation is particularly sensitive to pervasive and dense structural signals  that are often missed by $\mathcal{L}^\infty$-based procedures. 
   Moreover, we emphasize that our nonasymptotic approximation results hold for any finite sample size $n$, and they can also be applied to triangular arrays of random vectors $X_{1,n},\ldots, X_{n,n}$. 

For    statistical inference, particularly the feasible change-point detection procedure discussed in Section \ref{sec:changepoint}, a consistent estimator for the high-dimensional covariance matrix $\bm{\Sigma}$ is indispensable. We consider the following plug-in empirical estimator based on the Jackknife pseudo-values of the first-order projection. Specifically, we define the empirical projection for the $i$-th observation as:
\begin{align*}
    \widehat{\bm{g}}_i = \frac{1}{n-1} \sum_{j \neq i} \bm{h}(X_i, X_j) - \bm{U}_n, \quad \text{for } 1 \le i \le n. 
\end{align*}
The covariance matrix estimator is then constructed as:
\begin{align*}
    \widehat{\bm{\Sigma}} = \frac{1}{n} \sum_{i=1}^n \widehat{\bm{g}}_i \widehat{\bm{g}}_i^\top. 
\end{align*}
The following assumption and proposition establishes the consistency of $\widehat{\bm{\Sigma}}$. 
\begin{itemize}
    \item[(A3$^{\prime}$)] The bound $B$ in (A1) satisfies $B\asymp \sqrt{d}$ and the dimension $d=\o(n^{\frac{q-2}{2q}-\gamma})$ for some small $\gamma>0$.
\end{itemize}

\begin{proposition}\label{prop:sigma-consistency}
Suppose Assumptions (A1)-(A2) hold. Let $\| \cdot \|_{\text{op}}$ denote the operator norm, defined as $\|\mathbf{A}\|_{\text{op}} = \sup_{\|\bm{x}\|_2 = 1} \|\mathbf{A}\bm{x}\|_2$ for any matrix $\mathbf{A}$. Then the covariance matrix estimator $\widehat{\bm{\Sigma}}$ satisfies:
\begin{align*}
    \left\|\widehat{\bm{\Sigma}} - \bm{\Sigma}\right\|_{\text{op}} = \mathcal{O}_p\left( B^2 n^{ {1}/{q}- {1}/{2}} (\log n) \sqrt{\log d} + B^2 n^{ {2}/{q}-1} (\log n)^2 \log d + B \sqrt{ {d}/{n}} +  {d}/{n} \right). 
\end{align*}
Under Assumption (A3$^{\prime}$), the approximation error vanishes, yielding $\left\|\widehat{\bm{\Sigma}} - \bm{\Sigma}\right\|_{\text{op}} = \o_p(1)$ as $n \to \infty$.
\end{proposition}

It is worth emphasizing that the following lemma plays a key role in the proof of Theorem \ref{thm:sequential-borel}. By carefully exploiting the martingale structure of U-statistics and applying the multivariate martingale difference maximal inequality from \cite{Baiinequality}, we derive a sharp uniform maximal bound for the sequential partial sums of completely degenerate kernels. 

\begin{lemma}[Maximal Inequality for Degenerate U-statistics] \label{lem:degenerate-max}
    Let $X_1, \dots, X_n$ be independent random vectors taking values in $\R^p$. Let $\{\bm{f}_{ij}: \R^p \times \R^p \to \R^d\}_{1\le i<j\le n}$ be a sequence of measurable symmetric functions that are completely degenerate, i.e., $\E\{\bm{f}_{ij}(X_i, X_j) \mid X_i\} = \mathbf{0}$ and $\E\{\bm{f}_{ij}(X_i, X_j) \mid X_j\} = \mathbf{0}$ almost surely for all $i \neq j$. Assume (A2*) holds and define the degenerate partial sum process $\bm{M}_k = \sum_{1\le i<j\le k} \bm{f}_{ij}(X_i, X_j)$ for $2\le k\le n$. Then there exists a universal constant $C > 0$ such that:
    $$
    \E\left\{ \max_{2\le k\le n} \frac{\|\bm{M}_k\|_2}{(k-1)}  \right\} \le C D_2^2 d\log n.
    $$
\end{lemma} 
This bound  establishes the order of the remainder terms without invoking higher-order moment assumptions, and serves as a crucial technical building block for establishing the main theorems in this paper. 

\subsection{Gaussian approximation for full sums}
A byproduct of Theorem \ref{thm:sequential-borel} is introduced in this subsection, which is of independent interest. 
We now relax the identically distributed assumption to consider the independent but not necessarily identically distributed  scenario. Let $\{X_i\}_{i=1}^n$ be independent random vectors in $\R^p$ with respective distributions $\{F_i\}_{i=1}^n$. The symmetric kernel $\bm{h}(\cdot,\cdot)$ and the global U-statistic $\bm{U}_n$ are defined as before. However, the pairwise expectation now depends on the specific indices $i$ and $j$, i.e., $\bm{\theta}_{ij} := \E[\bm{h}(X_i,X_j)]$. We denote the generalized global parameter as $\bar{\bm{\theta}}_n := 2\sum_{1\le i<j\le n}\bm{\theta}_{ij}/{n(n-1)}.$ 
Then, the Hoeffding decomposition is accordingly modified to $$\bm{U}_n - \bar{\bm{\theta}}_n = \frac{2}{n}\sum_{i=1}^n \bm{g}_{i,n}(X_i) + \frac{1}{n(n-1)}\sum_{1\le i\neq j\le n} \bm{f}_{ij}(X_i,X_j)$$, where the projection term $\bm{g}_{i,n}(\cdot)$ and the degenerate kernel $\bm{f}_{ij}(\cdot,\cdot)$ are   given by:
\begin{align*}
\bm{g}_{i,n}(X_i) &= \frac{1}{n-1}\sum_{j \neq i}\left[\E\left\{\bm{h}(X_i,X_j)|X_i\right\}-\bm{\theta}_{ij}\right],\\
\bm{f}_{ij}(X_i,X_j) &= \bm{h}(X_i,X_j) - \E\left\{\bm{h}(X_i,X_j)|X_i\right\} - \E\left\{\bm{h}(X_i,X_j)|X_j\right\} + \bm{\theta}_{ij}.
\end{align*}
To establish the Gaussian approximation within this heterogeneous framework, we naturally adapt our regularity conditions as follows:
\begin{itemize}
\item[(A1*)] There exists some $q>2$ such that for each $1\le i\le n$, $\left(\E\|\bm{g}_{i,n}(X_i)\|^q\right)^{1/q}\le b_i$.
\item[(A2*)] Let $D_r = \max_{1\le i < j \le n}\max_{1\le m \le d} (\E|f_{m,ij}(X_i, X_j)|^r)^{1/r}$, and assume that $D_2 < \infty$.
\end{itemize}
\begin{theorem}\label{thm:fullsum-borel}
    Under Assumptions (A1*)-(A2*), for the standardized global statistic $\bm{T}_n := \sqrt{n}(\bm{U}_n - \bar{\bm{\theta}}_n)/2$, there exists a sequence of Gaussian vectors $\bm{W}_n =  \sum_{i=1}^n \bm Z_i/\sqrt{n}$, where $\bm{Z}_i \sim \mathcal{N}(\bm{0},\Cov\{\bm{g}_{i,n}(X_i)\})$, such that 
    $$    
    \left( \mathbb{E} \left\|\bm{T}_n-\bm{W}_n\right\|_2^2 \right)^{1/2} = \mathcal{O}\left((d/n)^{  {(q-2)}/{(2q)}}  \sqrt{\frac{\log n}{n} \sum_{i=1}^n b_i^2}\right).
    $$
\end{theorem}
Notably, the approximation error in Theorem \ref{thm:fullsum-borel} is governed by the  average of the projection moments, $n^{-1}\sum_{i=1}^n b_i^2$, rather than their maximum $\max_{1 \le i \le n} b_i^2$. This ensures the asymptotic validity of the strong Gaussian approximation even when a fraction of the observations come from highly heteroscedastic distributions. 


\section{Statistical applications}\label{Sec:statisticsapplication}
\subsection{Application to relevant tests}\label{sec:relevant}
In many statistical applications, rather than testing for exact equality, it is more practically meaningful to evaluate whether the structural distance between the parameters of two populations, $\bm{\theta}_1$ and $\bm{\theta}_2$, deviates by more than a scientifically acceptable tolerance level $\Delta > 0$ \citep{berger1987testing}, as further discussed and explored in \cite{dette2016detecting}. 
 The two-sample relevant hypotheses test is formulated as follows:
\begin{align*}
    H_0:\|\bm{\theta}_1 - \bm{\theta}_2\|_2^2 \le \Delta,\quad H_1:\|\bm{\theta}_1 - \bm{\theta}_2\|_2^2 > \Delta.
\end{align*}
This framework also covers the one-sample setting by treating one parameter as a pre-specified benchmark $\bm{\theta}_0$.


Consider that $\{X_i\}_{i=1}^{n_1}$ and $\{Y_i\}_{i=1}^{n_2}$ are two independent sequences of observations with sizes $n_1$ and $n_2$, respectively. Assume that $\min(n_1,n_2)\to\infty$ and $n_1/(n_1+n_2)\to\kappa\in(0,1)$. Let 
$N= {n_1n_2}/{(n_1+n_2)}$ 
denote the effective sample size, $\bm U_{k_1}^{(1)}$ and $\bm U_{k_2}^{(2)}$ denote the sequential U-statistics constructed from the first $k_1$ and $k_2$ observations of the two sequences, respectively. To construct a pivotal test without estimating the high-dimensional covariance matrix, we employ a self-normalization(SN) approach. Define the sequential distance process over the proportional time $\lambda \in (0, 1]$ as $$D_{n_1, n_2}(\lambda) := \left\|\bm{U}_{\lfloor \lambda n_1 \rfloor}^{(1)} - \bm{U}_{\lfloor \lambda n_2 \rfloor}^{(2)}\right\|_2^2.$$  For mathematical completeness, we set $D_{n_1, n_2}(\lambda) = 0$ whenever $\min(\lfloor \lambda n_1 \rfloor, \lfloor \lambda n_2 \rfloor) < 2$,  which includes the boundary $\lambda = 0$. The self-normalizer is then constructed from the centered partial-sum fluctuations of this process as
\[V_{n_1, n_2}^2 = \int_{0}^1 \left[ \lambda N \left\{D_{n_1, n_2}(\lambda) -   D_{n_1, n_2}(1) \right\}\right]^2 d\lambda.\]
The resulting SN test statistic is given by
\[S_{n_1, n_2} = \frac{N(D_{n_1, n_2}(1) - \Delta)}{V_{n_1, n_2}}.\]

To establish the asymptotic theory for $S_{n_1, n_2}$, let $\bm{\Sigma}_1$ and $\bm{\Sigma}_2$ be the covariance matrices of the first-order projections for $\bm h(X_i,X_j)$
and $\bm h(Y_i,Y_j)$, 
respectively. Write $\bm{\Sigma}_{\kappa} =  ({n_2}\bm{\Sigma}_1 +  {n_1} \bm{\Sigma}_2)/{(n_1+n_2)}$ and  $\bm{\delta} = \bm{\theta}_1 - \bm{\theta}_2$. The following regularity conditions are imposed for theoretical development.
\begin{itemize}
    \item [(B1)] The combined covariance matrix satisfies $\Tr(\bm{\Sigma}_{\kappa}) = \o\left(N^{1/2}(\log N)^{-1}\right)$.
    \item [(B2)] There exists a universal constant $c > 0$ such that $\bm{\delta}^\top \bm{\Sigma}_{\kappa} \bm{\delta} \ge c\|\bm{\delta}\|_2^2 > 0$ ($\bm\delta\neq \bm 0$). 
\end{itemize}

Assumption (B1) mildly constrains the variance of two U-statistics. 
Note that $\Tr(\bm{\Sigma}_{\kappa}) = \sum_{j=1}^d \lambda_j(\bm{\Sigma}_{\kappa})$, where $\lambda_j(\bm{\Sigma}_{\kappa})$'s are the eigenvalues of $\bm{\Sigma}_{\kappa}$. If all the eigenvalues are bounded, then (B1) effectively reduces to a requirement on the dimensionality, specifically $d = \o\left(N^{1/2} (\log N)^{-1}\right)$. Besides, Assumption (B2) guarantees that the quadratic form does not degenerate, and remains bounded away from zero whenever $\bm\delta\neq \bm 0$. This condition is a relaxation of the standard minimum eigenvalue assumption, i.e., $\lambda_{\min}(\bm{\Sigma}_{\kappa}) \ge c > 0$.

\begin{theorem}
\label{thm:main_sn}
Suppose Assumptions (A1)-(A3), and (B1)-(B2) hold for both independent samples. On the boundary of the null hypothesis that $ \|\bm{\delta}\|_2^2 = \Delta$, the SN test statistic is asymptotically pivotal:
\[S_{n_1, n_2} \xrightarrow{d} \mathbb{W} := \frac{\mathbb{B}(1)}{\left[ \int_0^1 \{\mathbb{B}(\lambda) - \lambda \mathbb{B}(1)\}^2 d\lambda \right]^{1/2}},\]
where $\mathbb{B}(\cdot)$ is a standard Brownian motion on $[0,1]$.
\end{theorem}

To establish the asymptotic power of the SN test, we consider a fixed alternative hypothesis where the true structural distance strictly exceeds the scientific tolerance threshold, i.e., $\Delta^*:=\|\bm{\delta}\|_2^2$, where $\Delta^* > \Delta$. 
The following Theorem \ref{thm:consistency} establishes the asymptotic consistency of the SN test.
\begin{theorem}
\label{thm:consistency}
  Under the alternative hypothesis $H_1: \|\bm{\delta}\|_2^2 = \Delta^* > \Delta$ and the   conditions in Theorem \ref{thm:main_sn}, if the following condition holds:
\[ \frac{\sqrt{N}(\Delta^* - \Delta)}{\sqrt{\bm{\delta}^\top \bm{\Sigma}_{\kappa} \bm{\delta}}} \to \infty \quad \text{as } \min(n_1, n_2) \to \infty, \]
then, 
 $S_{n_1, n_2} \xrightarrow{p} +\infty$. 
Furthermore, denote by $q_{1-\alpha}$ the $\alpha$-upper quantile of $\mathbb{W}$,  
\[\lim_{\min(n_1, n_2) \to \infty} \mathbb{P}(S_{n_1, n_2} > q_{1-\alpha}) = 1.\]
\end{theorem}

\bgroup
\renewcommand{\theexample}{\ref{exa:GMD}}
\begin{example}[Multivariate Gini's Mean Difference (continued)]
As a concrete application, we implement our self-normalized relevant testing procedure for the multivariate GMD to test $H_0: \|\bm{\theta}_1 - \bm{\theta}_2\|_2^2 \le \Delta$ versus $H_1: \|\bm{\theta}_1 - \bm{\theta}_2\|_2^2 > \Delta$. 
To see how our regularity conditions are mathematically evaluated in practice, consider the balanced sample case where the observations follow a multivariate Gaussian distribution, i.e., $X \sim \mathcal{N}(\bm{\mu}_X, \bm{\Sigma}_X)$ and $Y \sim \mathcal{N}(\bm{\mu}_Y, \bm{\Sigma}_Y)$. Let $\sigma_{X,\ell}^2 = (\bm{\Sigma}_X)_{\ell\ell}$ denote the marginal variances. The target parameter evaluates explicitly to $\theta_{1\ell} = \mathbb{E}|X_{1\ell} - X_{2\ell}| = 2\sigma_{X,\ell}/\sqrt{\pi}$.  
Under this specific distributional setting, the assumptions in our framework can be rigorously translated into structural constraints on the data generating process: 
\begin{enumerate}
    \item [\textbf{Assumption (A2):}] The second moment of the base kernel is analytically given by $\mathbb{E}|X_{1\ell} - X_{2\ell}|^2 = 2\sigma_{X,\ell}^2$. 
    Thus, Assumption (A2) holds strictly provided that the maximum marginal variance $\max_{1\le \ell \le d} \sigma_{X,\ell}^2 < \infty$.
    \item [\textbf{Assumption (B1):}] The combined covariance matrix $\bm{\Sigma}_{\kappa}$ is governed by the projection variances of both samples. For the GMD under Gaussianity, one can compute the projection variance via bivariate Gaussian integration, yielding $\Var(g_{1\ell}^{(X)}(X_{1\ell})) = c_g \sigma_{X,\ell}^2$ and $\Var(g_{1\ell}^{(Y)}(Y_{1\ell})) = c_g \sigma_{Y,\ell}^2$, where the absolute constant $c_g = 1/3 + (2\sqrt{3}-4)/\pi \approx 0.163$. For balanced samples, this establishes an exact linear scaling $\Tr(\bm{\Sigma}_{\kappa}) = {c_g} \left\{ \Tr(\bm{\Sigma}_X) + \Tr(\bm{\Sigma}_Y) \right\}/2$. Hence, Assumption (B1) dictates a bound on the underlying covariances: $\max\{ \Tr(\bm{\Sigma}_X), \Tr(\bm{\Sigma}_Y) \} = o(n^{1/2}(\log n)^{-1})$. If the marginal variances of both populations are bounded, this elegantly reduces to the dimension constraint $d = \o(n^{1/2}(\log n)^{-1})$.
\end{enumerate}
\end{example}
\addtocounter{example}{-1}
\egroup

\subsection {Application to   change point analysis}\label{sec:changepoint}

Beyond the relevant test discussed in Section \ref{sec:relevant}, our sequential Gaussian approximation framework naturally applies to change-point analysis. Formally, we aim to test the following hypothesis regarding the underlying parameter sequence $\bm{\theta}_t = \mathbb{E}[\bm{h}(X_t, X_{t'})]$, where $X_t, X_{t'}$ are independent copies at time $t$:
\begin{description}
    \item[$H_0$:]   $\bm{\theta}_1 = \bm{\theta}_2 = \cdots = \bm{\theta}_n = \bm{\theta}_0$.
    \item[$H_1$:] There exists an unknown change-point fraction $\tau^* \in (0, 1)$ with $k^* = \lfloor n \tau^* \rfloor$, such that $\bm{\theta}_1 = \cdots = \bm{\theta}_{k^*} \neq \bm{\theta}_{k^*+1} = \cdots = \bm{\theta}_n$.
\end{description}

To test whether the generalized moment $\bm{\theta}$ undergoes a structural change within the observation sequence, we construct a CUSUM-type statistic based on the sequential U-statistics. Let $k$ be the time index such that $2\le k\le n-2$. Recall that the sequential U-statistic $\bm{U}_k$ serves as an unbiased estimator for the target parameter using the pre-break subsample $\{X_1,\dots,X_k\}$. For the remaining post-break sample, we define 
$$\bm{U}_k^* = \binom{n-k}{2}^{-1} \sum_{k+1 \le i < j \le n} \bm{h}(X_i, X_j).$$
In this framework, we  adopt the ``first-versus-last'' CUSUM construction, directly contrasting the pre-break statistic $\bm{U}_k$ with the post-break statistic $\bm{U}_k^*$. Another alternative is the ``first-versus-full'' approach, which compares $\bm{U}_k$ against the global estimator $\bm{U}_n$. A comprehensive theoretical comparison by \cite{dehling2026} shows that, in the univariate parameter setting, these two paradigms are asymptotically equivalent under both the null hypothesis and sequences of local alternatives. Under fixed alternatives, their relative power depends  on the sign of the kernel's {\color{black}eccentricity}, implying that neither approach universally dominates the other.


Formally, the sequential CUSUM process $\bm{\mathcal{C}}_n(k)$ is defined as the scaled difference between $\bm U_k$ and $\bm U_k^\ast$:  
\begin{align*}
    \bm{\mathcal{C}}_n(k) &:= \sqrt{n} \frac{k}{n} \left( 1 - \frac{k}{n} \right) (\bm{U}_k - \bm{U}_k^*)\\
    &= \bm{\mathcal{D}}_n(k) + \frac{1}{\sqrt{n}} \left\{ \sum_{i=1}^{k} 2\bm{g}(X_i) - \frac{k}{n} \sum_{i=1}^n 2\bm{g}(X_i) \right\} + \bm{\Delta}_n(k).
\end{align*} 
Here, the deterministic drift term $\bm{\mathcal{D}}_n(k)$ captures the structural shift in the underlying distribution and is defined as: 
$$
\bm{\mathcal{D}}_n(k) = \sqrt{n} \frac{k}{n} \left( 1 - \frac{k}{n} \right) \left\{ \mathbb{E}(\bm{U}_k) - \mathbb{E}(\bm{U}_k^*) \right\}.
$$
Under the null hypothesis $H_0$, the expectations coincide, yielding $\bm{\mathcal{D}}_n(k) \equiv \bm{0}$. Under the alternative hypothesis $H_1$, this drift term provides the signal that drives the test statistic to diverge. 
The function $\bm{g}(X_i)$ consistently represents the centered first-order projection of the U-statistic kernel as before and the remainder process $\bm{\Delta}_n(k)$ captures the stochastic contribution from the completely degenerate kernel components $\bm{f}(X_i, X_j)$, and is given by:
$$
\bm{\Delta}_n(k) = \frac{2\sqrt{n}}{n^2} \left\{ \frac{n-k}{k-1} \sum_{1 \le i < j \le k} \bm{f}(X_i, X_j) - \frac{k}{n-k-1} \sum_{k+1 \le i < j \le n} \bm{f}(X_i, X_j) \right\}.
$$

\begin{remark}[{\color{black}Compared with \cite{wegner2023robust}}]
    It is worth noting the distinction between our statistic $\bm{\mathcal{C}}_n(k)$ and the U-statistic-based CUSUM process $\mathcal{W}_{n,k} = n^{-3/2} \sum_{i=1}^k \sum_{j=k+1}^n \bm{h}^{anti}(X_i, X_j)$ proposed by \cite{wegner2023robust}. Notably, $\mathcal{W}_{n,k}$ is a two-sample U-statistics utilizing an antisymmetric kernel, which takes the first $k$ and last $n-k$ observations as two population. This design   targets the overall distributional equivalence by testing $\mathbb{E}[\bm{h}^{anti}(X,Y)]=\bm{0}$. Conversely, $\bm{\mathcal{C}}_n(k)$ contrasts two separate one-sample U-statistics based on a symmetric kernel,   monitoring a target structural parameter $\bm\theta = \mathbb{E}[\bm{h}(X,X')]$. 
    Despite these structural and conceptual differences, their asymptotic mechanisms are equivalent. By the Hoeffding decomposition, both processes are dominated by analogous linear bridges driven by their first-order projections. Moreover, our symmetric framework strictly encompasses the classical observation-driven CUSUM process when specifying the linear kernel $h(x,y) = x+y$.
\end{remark}

To test $H_0$, we   define the $\mathcal{L}^2$-type test statistic as $T_n = \max_{2 \leq k \leq n-2} \|\bm{\mathcal{C}}_n(k)\|_2$. Deriving its asymptotic null distribution relies  on the theoretical tools provided in Section \ref{sec:sequential}. By employing a martingale-based bound for the degenerate remainder $\bm{\Delta}_n(k)$, 
we impose the following mild trace bound:
\begin{itemize}
\item [(C1)] The covariance matrix satisfies $\Tr(\bm{\Sigma}) = \o\left(n(\log n)^{-1}\right)$, where $\bm{\Sigma}$ is defined as before.
\end{itemize}
The following theorem establishes the uniform Gaussian coupling for the proposed sequential process $\bm C_n(\cdot)$, which   yields the limiting behavior of $T_n$.\begin{theorem}\label{thm:cusum_ga}Suppose Assumptions (A1)-(A3) and (C1) hold. Under $H_0$, on a sufficiently rich probability space, there exists a sequence of Gaussian processes $\{\bm{\mathcal{B}}_n(t)\}_{t \in (0,1)}$, where each $\bm{\mathcal{B}}_n$ is distributed as a $d$-variate Brownian bridge with covariance function $\mathbb{E}[\bm{\mathcal{B}}_n(s)\bm{\mathcal{B}}_n(t)^\top] = 4(\min(s,t) - st)\bm{\Sigma}$, such that:$$\sup_{t \in (0, 1)} \left\| \bm{\mathcal{C}}_n(\lfloor nt \rfloor) - \bm{\mathcal{B}}_n(t) \right\|_2 = \o_p(1).$$
Hence, $$\left| T_n - \sup_{t \in (0, 1)} \| \bm{\mathcal{B}}_n(t) \|_2 \right| = \o_p(1), \quad \text{as } n \to \infty.$$\end{theorem}

While Theorem \ref{thm:cusum_ga} establishes the asymptotic null distribution of $T_n$, its dependence on the unknown $\bm{\Sigma}$ makes the direct use of $\bm{\mathcal{B}}_n(t)$ infeasible. To construct a practical testing procedure, we substitute $\bm{\Sigma}$ with the Jackknife estimator $\widehat{\bm{\Sigma}}$ established in Section \ref{sec:sequential}. However, ensuring the uniform validity of this substitution requires a slightly stronger restriction on the dimension growth rate:
\begin{itemize}
    \item[(A3$^{\prime\prime}$)] The bound $B$ in (A1) satisfies $B\asymp \sqrt{d}$ and the dimension $d=\o\left(n^{\frac{q-2}{3q}-\gamma}\right)$ for some small $\gamma>0$.
\end{itemize}
\begin{corollary}\label{cor:BB-coupling}
Assume that the minimal eigenvalue of the  covariance matrix is strictly bounded away from zero, i.e., $\lambda_{\min}(\bm{\Sigma}) \ge c > 0$ for some constant $c > 0$. Let $\widehat{\bm{\mathcal{B}}}_n(t)$ denote the d-dimensional Brownian bridge which has the covariance function $\mathbb{E}[\widehat{\bm{\mathcal{B}}}_n(s)\widehat{\bm{\mathcal{B}}}_n(t)^\top] = 4(\min(s,t) - st)\widehat{\bm{\Sigma}}$. Under Assumptions (A1)-(A2) and (A3$^{\prime\prime}$), 
\begin{align*}
    \sup_{t \in (0,1)} \left\| \widehat{\bm{\mathcal{B}}}_n(t) - \bm{\mathcal{B}}_n(t) \right\|_2 = \o_p(1).
\end{align*}
\end{corollary}

Building upon the uniform coupling in Corollary \ref{cor:BB-coupling}, we propose a resampling  procedure to empirically  compute the  critical values. Conditional on the observed data, we generate a large number of independent paths of a standard $d$-dimensional Brownian bridge $\mathcal{B}(t)$. For each path, we computationally realize the feasible Gaussian process via the linear transformation $\widehat{\bm{\mathcal{B}}}_n(t) = 2 \widehat{\bm{\Sigma}}^{1/2} \mathcal{B}(t)$. The feasible critical value $\hat{q}_{1-\alpha, n}$ is then obtained as the upper-$\alpha$ quantile of the simulated distribution of $\sup_{t \in (0,1)} \left\|\widehat{\bm{\mathcal{B}}}_n(t)\right\|_2$.
Therefore,  rejecting the null hypothesis when $T_n > \hat{q}_{1-\alpha, n}$ guarantees an asymptotic size exactly equal to $\alpha$ under $H_0$.

We now investigate the asymptotic behavior of the test statistic and the change-point estimator under the alternative hypothesis $H_1$. 
Consider that $\{X_i\}_{i=1}^{k^*}$ are i.i.d. following a distribution $F_1$, and $\{X_i\}_{i=k^*+1}^{n}$ are i.i.d. following another distribution $F_2$, leading to a parameter shift:
$$
\bm{\theta}_t = \begin{cases} \bm{\mu}_1, & 1 \le t \le k^*, \\ \bm{\mu}_2, & k^* < t \le n, \end{cases}
$$
with a  non-zero shift size $\bm{\delta} = \bm{\mu}_2 - \bm{\mu}_1 \neq \bm{0}$. Let $\bm{\Sigma}_1$ and $\bm{\Sigma}_2$ denote the covariance matrices of the centered first-order projections under the respective distributions $F_1$ and $F_2$, defined analogously to $\bm{\Sigma}$. We define $\Tr_{\max} := \max\{\Tr(\bm{\Sigma}_1), \Tr(\bm{\Sigma}_2)\}$. The following theorem establishes the asymptotic consistency of the proposed test.
\begin{theorem}\label{thm:test_consistency}
  Assume that  $\Tr_{\max} = \o\left(n(\log n)^{-1}\right)$. Under the alternative hypothesis $H_1$ and Assumptions (A1)-(A3), if the shift magnitude satisfies $\sqrt{n}\|\bm{\delta}\|_2 / \sqrt{\Tr_{\max}} \to \infty$ as $n \to \infty$, then, the test has asymptotic power 1, 
$$\lim_{n \to \infty} \mathbb{P}(T_n > \hat{q}_{1-\alpha, n}) = 1.$$
\end{theorem}

Beyond  detecting the presence of a structural break, it is often of   interest to precisely locate it. We naturally define the change-point estimator as the time index maximizing the Euclidean norm of the proposed CUSUM process:
$$
\hat{k} = \operatorname*{arg\,max}_{2 \le k \le n-2} \|\bm{\mathcal{C}}_n(k)\|_2.
$$
The following theorem guarantees that the estimator $\widehat k$ consistently locates the true change-point fraction.

\begin{theorem}
\label{thm:estimator_consistency}
Suppose the conditions of Theorem \ref{thm:test_consistency} hold. Let $\bm{\mu}_{12} = \mathbb{E}[h(X, Y)]$ denote the cross-expectation where $X \sim F_1$ and $Y \sim F_2$. Assume $\bm{\mu}_{12}$ satisfies the following geometric constraints:
\begin{equation*}
\left\{
\begin{aligned}
    &2\tau^* (\bm{\mu}_2 - \bm{\mu}_1)^\top (\bm{\mu}_1 - \bm{\mu}_{12}) < \|\bm{\delta}\|_2^2, \\
    &2(1-\tau^*) (\bm{\mu}_2 - \bm{\mu}_1)^\top (\bm{\mu}_{12} - \bm{\mu}_2) < \|\bm{\delta}\|_2^2.
\end{aligned}
\right.
\end{equation*}
Then, the change-point estimator $\hat{k}$ is consistent for the true change-point fraction $\tau^*$. That is,
$$
\hat{\tau} = {\hat{k}}/{n} \xrightarrow{\P} \tau^* \quad \text{as } n \to \infty.
$$
\end{theorem}

The geometric constraints introduced in Theorem \ref{thm:estimator_consistency} are  required to ensure that the asymptotic population drift function $V(t)$ possesses a unique global maximum exactly at the true change-point $\tau^*$. As   shown in the proof, the rescaled objective function uniformly converges to the deterministic limit $V(t)$, which  takes the following form:
\begin{equation*}
    V(t) = \begin{cases} 
    \left\| \frac{t}{1-t} \left[ (1-\tau^*)^2 \frac{-\bm{\delta}}{\|\bm{\delta}\|_2} + 2(\tau^*-t)(1-\tau^*) \frac{\bm{\mu}_1 - \bm{\mu}_{12}}{\|\bm{\delta}\|_2} \right] \right\|_2, & t \le \tau^*, \\
    \left\| \frac{1-t}{t} \left[ (\tau^*)^2 \frac{-\bm{\delta}}{\|\bm{\delta}\|_2} + 2\tau^*(t-\tau^*) \frac{\bm{\mu}_{12} - \bm{\mu}_2}{\|\bm{\delta}\|_2} \right] \right\|_2, & t > \tau^*.
    \end{cases}
\end{equation*}
In fact, these specific geometric inequalities can be relaxed: any milder condition guaranteeing that $V(t)$ attains a strict and unique supremum at $t = \tau^*$ is  sufficient for the consistency result to hold. For analogous geometric conditions and an extensive discussion in the univariate setting, we refer to \cite{dehling2026}.

\bgroup
\renewcommand{\theexample}{\ref{exa:CDP}}
\begin{example}[Characteristic Dispersion Parameter (continued)]
To illustrate the practical relevance of our theoretical results, we revisit the characteristic dispersion parameter $\bm{h}(X_i, X_j) = (\cos(X_{i1} - X_{j1}), \dots, \cos(X_{id} - X_{jd}))^\top$. 
In quantitative finance, detecting structural shifts in market volatility is crucial. However, financial returns often exhibit heavy tails and extreme outliers (\cite{econbg}), causing classical variance-based change-point methods to produce frequent false alarms. Our CUSUM procedure relieves this by targeting the robust dispersion parameter $\bm{\theta} = (\theta_1, \dots, \theta_d)^\top$, where $\theta_\ell = \E[\cos(X_{i\ell} - X_{j\ell})]$, which serves as a stable proxy for generalized volatility.

Theoretically, this robustness is guaranteed by the geometric boundedness of this kernel. As shown via the Hoeffding decomposition, the first-order projection and the degenerate kernel are universally bounded. This ensures that the maximal norm conditions in Assumptions (A1)-(A2) are satisfied  without imposing any other moment constraints on the underlying distribution.  
By tracking these strictly bounded components, our framework   filters out the noise of unbounded market shocks, offering a robust tool for econometric structural break analysis.
\end{example}
\addtocounter{example}{-1}
\egroup

\bgroup
\renewcommand{\theexample}{\ref{exa:SKT}}
\begin{example}[Spatial Kendall’s Tau Matrix (continued)]
We now deploy the spatial Kendall's tau matrix introduced in Example \ref{exa:SKT} within our sequential testing framework to monitor the structural evolution of gene regulatory networks. 
In single-cell studies, identifying critical transitions, such as cell-fate decisions along a pseudotime trajectory, requires detecting structural rewiring in the underlying gene co-expression network. Classical covariance-based change-point tests often produce spurious breaks due to temporal transcriptional bursting. Our CUSUM procedure naturally overcomes this by monitoring the robust spatial association. 
Since $\|\bm{h}(X_i, X_j)\|_2 = 1$,  the maximal norm conditions   in Assumptions (A1)-(A2) are  again satisfied without imposing any  moment constraints.  
Therefore, our CUSUM statistic operates as a purely scale-invariant change-point detector. It ensures that any detected structural break represents a genuine topological shift in the genetic dependence structure and is robust to large single-cell measurement outliers. 
\end{example}
\addtocounter{example}{-1}
\egroup


\section{Concluding remark}\label{Sec:conclude}

This paper develops a strong Gaussian approximation theory for high-dimensional, vector-valued $U$-statistics of order two. The main probabilistic result is a sequential coupling in Euclidean norm between the centered and scaled $U$-statistic process and a Gaussian partial-sum process, with an explicit error bound that vanishes under polynomial growth of the dimension. The proof combines the Hoeffding decomposition with a high-dimensional strong approximation for independent sums and, crucially, a maximal inequality for vector-valued completely degenerate $U$-statistics obtained through a martingale argument. In this way, the linear projection and the degenerate remainder are treated within a unified framework, yielding a pathwise approximation for the whole sequential process rather than only a weak approximation for a fixed functional.

Beyond the main coupling theorem, the paper also provides several complementary results that strengthen the overall framework. We derive a Gaussian approximation for the global statistic under independent but not necessarily identically distributed observations, and we establish consistency of a jackknife-type estimator for the covariance matrix of the first-order projection. These probabilistic tools are then translated into concrete inferential procedures. On the one hand, we obtain Brownian-bridge approximations for $U$-statistic based CUSUM processes, together with feasible critical values and consistency of the corresponding change-point method. On the other hand, we develop a self-normalized relevant test for hypotheses of the form $\|\bm{\theta}-\bm{\theta}_0\|_2^2\le \Delta$, whose null limit is pivotal and therefore avoids estimating a high-dimensional covariance matrix directly. A notable advantage of the present framework is that it naturally accommodates bounded kernels, so that the theory continues to apply under heavy-tailed sampling where conventional moment-based procedures may become unstable or even ill posed.

At the same time, several limitations of the present work should be noticed. First, the main sequential theory is developed for non-degenerate $U$-statistics of order two under independence. Although this setting already covers a broad range of robust kernels and important applications, it does not include genuinely degenerate problems or more general dependence structures. Second, the dimension regime allowed by our strong approximation is polynomial rather than exponential. This reflects the fact that we work in the $\mathcal{L}^2$-geometry and seek a uniform-in-time coupling, which is well suited to dense alternatives and energy-type functionals but is not designed for extremely sparse signals in ultra-high dimensions, where $\mathcal{L}^\infty$-type Gaussian approximations are often more effective. Third, the feasible change-point procedure still relies on covariance estimation together with trace and eigenvalue conditions, and these assumptions may become restrictive when the covariance structure is nearly singular or when the dimension is extremely large relative to the sample size.

One important direction for future research is to extend the present theory to dependent and locally nonstationary observations. From the viewpoint of the linear projection, such an extension appears plausible in light of the existing strong approximation theory for high-dimensional partial sums. The real obstacle lies in the canonical $U$-statistic remainder, whose temporal dependence destroys the simple martingale structure exploited here. Developing maximal inequalities and strong couplings for degenerate $U$-statistics under mixing, physical dependence, or locally stationary regimes would substantially broaden the scope of the present results and would be particularly relevant for time-series and functional-data change-point problems.

A second important extension is to move beyond order-two kernels. Higher-order $U$-statistics, $V$-statistics, and incomplete or size-dependent kernels arise naturally in semiparametrics, network analysis, and computationally scalable high-dimensional procedures. In such settings, the Hoeffding decomposition contains several projection levels, so that a sequential strong approximation would require simultaneous control of multiple canonical components with explicit dimension dependence. Establishing such a theory would connect the present Euclidean strong-approximation viewpoint with a broader class of nonlinear high-dimensional statistics and could lead to new inferential tools beyond the two applications studied in this paper.

Overall, we hope that the present work serves as a first step toward a broader strong approximation theory for nonlinear statistics in high dimensions, in which coupling arguments play a central role not only in probability theory but also in the construction of robust and practically feasible inferential procedures.

\bibliographystyle{apalike}
\bibliography{ref_arxiv}

\end{document}